\documentclass[12pt]{amsart}
\usepackage{amssymb,amsfonts,amsmath}
\usepackage[colorlinks]{hyperref}

\newcommand{\br}[1]{\left[{#1}\right]} 
\newcommand{\brf}[1]{\left\{{#1}\right\}}
\newcommand{\brr}[1]{\left({#1}\right)} 

\newcommand{\setR}{{\mathord{\mathbb R}}} 

\newcommand{\va}{\mathbf{a}} 
\newcommand{\vb}{\mathbf{b}} 
\newcommand{\vc}{\mathbf{c}} 
\newcommand{\vd}{\mathbf{d}} 
\newcommand{\vx}{\mathbf{x}} 
\newcommand{\vy}{\mathbf{y}} 
\newcommand{\vz}{\mathbf{z}} 
\newcommand{\vu}{\mathbf{u}} 
\newcommand{\vv}{\mathbf{v}} 

\newcommand{{\0}}{\mathbf{0}}

\newtheorem{Thm}{Theorem}

\newcommand{\bproof}{\noindent{\bf Proof \hspace{1 mm}}}
\newcommand{\eproof}{\hfill$\Box$}

\title[Applications of Matrices
Multiplication]
{Applications of Matrices
Multiplication  to Determinant and Rotations formulas in
$\setR^n$}

\date{}

\author[A. Goldvard]{Alex Goldvard}
\email{goldvard@braude.ac.il}
\address{Department of Mathematics,
         ORT Braude College,
         P.O. Box 78,
         21982 Karmiel,
         Israel}
\author[L. Karp]{Lavi Karp}
\email{karp@braude.ac.il}
\address{Department of Mathematics,
         ORT Braude College,
         P.O. Box 78,
         21982 Karmiel,
         Israel}


\begin{document}
\maketitle
 \begin{abstract}
This note deals with two topics of linear algebra. We give a simple and short
proof of the multiplicative property of the determinant and provide a
constructive formula for  rotations. The derivation  of the rotation matrix
relies on simple matrix calculations and thus can be presented in an elementary
linear algebra course. We also classify all invariant subspaces of equiangular
rotations in 4D.
\end{abstract}

\section{Introduction}

This article aims to  promote  several geometric aspects of linear algebra. The
geometric motivation often leads   to simple proofs in addition to  increasing
the students' interest  of this subject and providing a solid basis. Students
with a confident grasp of these ideas will encounter little difficulties in
extending them  to more abstract linear spaces. Many geometrical operations can
be rephrased in the language of vectors and matrices. This includes projections,
reflections and rotations, operations which have numerous applications in
engineering, physics, chemistry and economic, and  therefore their matrix's
representation should be included in basic linear algebra course.

There are basically two attitudes in  teaching linear algebra. The abstract one
which deals  with formal definitions of vector spaces, linear transformations
etc. Contrary to the abstract vector spaces, the analytic approach   deals
mainly with the vector space $\mathbb{R}^n$ and provides the basic concepts and
proofs in these spaces.

However, when the  the analytic approach deals with the definition of linear
transformations it uses the notion of {\it  representation of the matrix of the
transformation} in arbitrary basis and thus it actually goes back to the
abstract setting. To clarify this issue, let us consider for example the
calculations a rotation $T$ of a vector $\vx$ in $\setR^3$. In this way one
start with choosing an appropriate basis $B$, and calculating the matrix of the
given transformation in this basis $[T]_B$. In the second step one has to
compute the transformation matrix   $P$ from the standard basis to the basis $B$
and its inverse $P^{-1}$. The final step consists of applying the matrix $
P^{-1}[T]_BP$ to $\vx$.  This cumbersome machinery is common for both the
abstract and the analytic approaches. In addition, it take an essential effort
to teach all the necessary details in order to use this non-practicable formula.

It is amazing why one should use this complicate formula while the Rodrigues'
rotation formula does  it  efficiently.  Although its simplicity, Rodrigues'
formula does not appear in the current linear algebra textbooks.  Its proof is
elementary, but require non-trivial geometric insight.

In this note  we present a new  proof of  the matrix's representation of
Rodrigues' formula. The essential point is that we regard the multiplication
$A\vx$, of a vector $\vx$ by a matrix $A$,  simultaneously as an algebraic
operation and geometrical transformation (similarly to Lay \cite{Lay}). This
enable us to derive the rotation formula in the three dimensional space. In
higher dimensional spaces we first propose a geometric definition of a rotation
and after that we derive the formula in a similar manner to the three
dimensional spaces.   Having calculate the matrix of rotation according to that
definition, we show it  is identical to the common definition of rotation, that
is, an orthogonal matrix with determinant one.

We also consider the multiplication of two matrices $AB$ as multiplications of
the columns of $B$ by the matrix $A$. {Applying this point of view we provide} a
simple proof of the multiplicative property of the determinants. The standard
proof of  this property is often being skipped from the class room since it is
considered as too complicated. The proof which we present here could easily be
thought in the beginning of a linear algebra course.

\section{Basic facts and notations }
\label{sec:notations}

We recall the definitions of multiplication of a vector by a matrix
and the multiplications of two matrices. {Both
definitions rely solely }  on the basic two operations
of vectors in  $\setR^n$, namely, addition and multiplication by a scalar.

 We denote $m\times n$ matrix $A$ by $\br{\va_1,\va_2,\ldots ,\va_n}$,
where $\{\va_1,\va_2,\ldots ,\va_n\}$ are the columns of $A$. Let
$\vx=\begin{bmatrix} x_1 \\ \vdots \\ x_n \end{bmatrix}$ be a vector in
$\setR^n$,
then
\begin{equation}
\label{MVM}
A\vx = x_1\va_1 + x_2\va_2 + \ldots + x_n\va_n.
\end{equation}
This means that $A\vx$ is simply a linear combination of
the columns of the matrix $A$. Conversely, any linear combination of $n$ vectors
$\{\va_1,\va_2,\ldots ,\va_n\}$ can be written as a matrix multiplication.
Note that beside of $A\vx$ being a linear combination,  we can interpret
it as a
transformation from $\setR^n$ to $\setR^m$
by corresponding to  each $\vx\in\setR^n$ the vectors $A\vx = \vy \in
\setR^m$. Obviously this operation has the
linearity property:
\begin{equation}
\label{Lin}
 A(\alpha\vu+\beta\vv)=\alpha A\vu + \beta A\vv.
\end{equation}
Thus any linear transformation
can be written as a multiplication of vectors by a matrix
and therefore the formal definition the linear transformations
{from $\setR^n$ to $\setR^m$ } seems to be  superfluous.

Most of the textbooks define a matrix multiplications by the row-column rule.
But
 the original definition of Cayley is by
means of a composition
of two linear substitutions (see e.g. \cite{Me,Tu}). This means that if
$A$
is a $m\times n$ matrix and $B$
$n\times k$, then the matrix $C=AB$ is defined through the identity
$C\vx =
A\brr{B\vx}$,
where  $\vx\in\setR^k$. From (\ref{MVM}) {and (\ref{Lin}) }it  immediately
follows
that
\begin{equation}\label{CD}
	AB = \br{A\vb_1,A\vb_2,\ldots , A\vb_k},
\end{equation}
where $\{\vb_1,\vb_2,\ldots ,\vb_k \}$ are the columns of the matrix $B$.

Another useful way to multiply matrices is by the column row rule, that is,
\begin{equation}\label{CRM}
	AB = \va_1\vb_1^T + \ldots + \va_n\vb_n^T,
\end{equation}
where $\vb_1^T,\ldots ,\vb_n^T$ are rows vectors of the transpose
$B^T$. This type of multiplication will be used in Section \ref{sec:rotation}.

\section{Multiplicative property of determinant}
\label{sec:det}

{Let $A$ be $n\times n$ matrix
with coefficients $[a_{ij}]_{ij=1}^n$.}  The determinant of $A$
is defined to be the scalar
\begin{equation}
\det(A) = \sum_p	\sigma(p)a_{1p_1}a_{2p_2}\ldots
a_{np_n},
\end{equation}
where the sum is taken over the $n!$ permutations $p=\brr{p_1,p_2,\ldots p_n}$
of $\brr{1,2,\ldots ,n}$ and
\[ \sigma(p)=\begin{cases}
	+1 \; \mathrm{if} \, p\, \text{is\, even permutation} \\
	-1 \; \mathrm{if} \,p\,  \text{is\, odd permutation}
\end{cases}. \]

\begin{Thm}
\label{MPD1}
Suppose $A$ and $\br{\vu,\vv,\ldots ,\vz}$ are two $n \times n$
matrices, then
\begin{equation}
\label{eq6}
\det \left(\br{A\vu,A\vv,\ldots ,A\vz}\right) = \det (A )
\det (\br{\vu,\vv,\ldots ,\vz}).
\end{equation}
\end{Thm}
\bproof Denote {$A$ by
$\br{\va_1,\va_2, \ldots ,\va_n}$} and  {let}
\[ \vu = \begin{bmatrix}
	u_1 \\ u_2 \\ \vdots \\ u_n
\end{bmatrix},
\vv = \begin{bmatrix}
	v_1 \\ v_2 \\ \vdots \\ v_n
\end{bmatrix},
\vz = \begin{bmatrix}
	z_1 \\ z_2 \\ \vdots \\ z_n
\end{bmatrix}.  \]
{We are now
using  formula (\ref{MVM}), the linearity of the determinant and the fact that a
matrix
having two equal columns its determinant is zero. All these result with}
\begin{gather*}
	\det (\br{A\vu,A\vv,\ldots ,A\vz}) = \det \left(\br{\sum_{i=1}^n u_i
\va_i ,
 \sum_{i=1}^n v_i \va_i ,\ldots , \sum_{i=1}^n z_i \va_i }\right) \\ =
	\sum_p (u_{p_1} v_{p_2} \ldots z_{p_n})\det( \br{\va_{p_1},\va_{p_2},
\ldots ,\va_{p_n}}).
\end{gather*}
Since the determinant changes sign when two columns are interchanged, we
have
\[ \det( \br{\va_{p_1},\va_{p_2},\ldots ,\va_{p_n}}) = \sigma(p) \det (A). \]
Therefore,
\[\begin{split}
&\det (\br{A\vu,A\vv,\ldots ,A\vz}) = \det (A) \sum_p \sigma(p)u_{p_1} v_{p_2}
\ldots z_{p_n} \\=& \det (A)  \det (\br{\vu,\vv,\ldots ,\vz}).
  \end{split}
\]
\eproof

The above Theorem has an important geometric interpretation
and we discuss it here
in $\setR^2$. It is well known that if $\vu,\vv$ are two collinear vectors, then
$\det (\br{\vu,\vv})$ is the area of parallelogram spanned by $\{\vu,\vv\}$. Now
if $\det (A) \neq 0$, then $\{A\vu, A\vv\}$  span also a parallelogram.
Therefore
the number $\det (A)$ is the proportion between the areas of 
 the parallelograms spanned  by $\{\vu,\vv\}$ and $\{A\vu,A\vv\}$.

Note that Theorem \ref{MPD1} gives  the multiplicative
property of the determinant. Indeed, let  $  B = \br{\vu,\vv,\ldots ,\vz}$,
then $ AB=\br{A\vu,A\vv,\ldots ,A\vz}  $ and (\ref{eq6}) becomes
\[  \det \brr{AB} = \det (A) \det (B). \]

\section{The derivation of a formula for \\ a rotation matrix in $\setR^n$}
\label{sec:rotation}

First of all let us discuss the definition of a rotation in $\setR^n$. The
common definition of a rotation is by means of an orthogonal matrix with
determinant one. Here we provide another definition based on geometric
considerations. We are aware  that such definition was probably given in the
past, but we could note traced it. Its advantage is being practicable and the
computation of the rotation's matrix does not  utilize   eigenvalues and
eigenvectors. We shall then verify the equivalence of the two definitions.

  It is rather
simple to define a rotation in $\setR^2.$ A linear transformation
$R$ is a rotation if the angle between the vectors $R\vx$ and $\vx$ is a
constant
for all $\vx\in\setR^2.$ From this definition follows that
\[ R\begin{bmatrix}
	1 & 0 \\
	0 & 1
\end{bmatrix} =
\begin{bmatrix}
	\cos\alpha & -\sin\alpha \\
	\sin\alpha  & \cos\alpha
\end{bmatrix}, \]
where the rotation is by an angle $\alpha$ counterclockwise. Therefore
\begin{equation}\label{R2}
R = \begin{bmatrix}
	\cos\alpha & -\sin\alpha \\
	\sin\alpha  & \cos\alpha
\end{bmatrix} = \cos \alpha I +\sin\alpha J
\end{equation}
{where $I$ is the identity matrix and} $J=\begin{bmatrix} 0 &-1\\ 1& 0
\end{bmatrix}$.
Formula (\ref{R2}) implies that $||R\vx|| = ||\vx||$ for each $\vx\in\setR^2$
 and that confirms our geometric intuition. Matrices which
preserve norm are called orthogonal matrices. Their
determinant is $\pm 1$ and rotations are orthogonal matrices with
determinant one.

The  definition of a rotation in $\setR^3$ is slightly more
involved. A linear transformation $R$ is called rotation
if there exists two dimensional subspace $\Pi$ of $\setR^3$ (the plane of the
rotation) such that  the angle between vectors $R\vx$ and
$\vx$ is a constant for all $\vx\in\Pi$, and $R\vy = \vy$ for each
$\vy$ orthogonal to $\Pi$ (the axis of
the rotation). Euler's theorem about the rigid motion of the sphere with fixed
center justifies this definition.

In order to calculate the matrix $R$ we pick two {orthonormal} vectors
$\va,\vb\in\Pi$ and a unit vector $\vc$ which is orthogonal to $\Pi$ and  such
that the triple
$\{\va,\vb,\vc\}$ is right-handed.  Applying the rotation $R$ to these vectors
results with
\begin{equation*}
R \va  =  \cos\alpha \va +\sin\alpha \vb, \ \
R \vb =  -\sin\alpha \va + \cos\alpha \vb,\ \
R \vc =  \vc .
\end{equation*}
We write the above equalities in a matrix form $RP=Q$, where
$P = \br{\va, \vb, \vc}$ and $ Q = \br{R \va, R \vb, R \vc} $.
Since $P$ is an orthogonal matrix, $R=QP^T$ and calculating $QP^T$
by means of (\ref{CRM}), we get that
\begin{equation}\label{eq2}
\begin{split}
R &=\br{R\va,R\vb,R\vc}\begin{bmatrix}
	\va^T \\
	\vb^T \\
	\vc^T
\end{bmatrix} = (R\va)\va^T + (R\vb)\vb^T + (R\vc)\vc^T \\
&=\cos\alpha\brr{\va\va^T+\vb\vb^T}+\sin\alpha\brr{\vb\va^T-\va\vb^T}
+\vc\vc^T.
\end{split}
\end{equation}

The skew symmetric matrix  $(\vb\va^T-\va\vb^T)$  is the matrix representation
of the cross product $\vc\times \vx$. To see this note that
$(\vb\va^T-\va\vb^T)\vc=\0$, $(\vb\va^T-\va\vb^T)\va=\vb$ and
$(\vb\va^T-\va\vb^T)\vb=-\va$. Hence
\begin{equation*}
 R\vx=\cos\alpha \vx+(1-\cos\alpha)\vc\vc^T\vx+\sin\alpha (\vc\times \vx),
\end{equation*}
which is the known Rodrigues' formula.

Formula (\ref{eq2}) resembles to a large extent the two dimensional formula
(\ref{R2}).
The matrix $(\va\va^T+\vb\vb^T)$ is the projection on the plane $\Pi$,
$\vc\vc^T$ is the projection on the line orthogonal to $\Pi$ and
\begin{equation}
\label{eq:9}
 (\vb\va^T-\va\vb^T)^2=-(\va\va^T+\vb\vb^T). 
\end{equation} 
Since the rotation is actually in the plane $\Pi$, we see that 
$(\va\va^T+\vb\vb^T)$ corresponds $I$ and $(\vb\va^T-\va\vb^T)$ corresponds $J$
in formula (\ref{R2}). In both formulas, $\cos\alpha$ is the coefficient of a
symmetric matrix and $\sin\alpha$  is the  coefficient of an
anti-symmetric matrix.

It is easy to check  that the matrix $R(\alpha):=R$
in (\ref{eq2}) is an
orthogonal matrix with determinate one. {Indeed,} relation (\ref{eq:9}) and the
orthogonality of  $\{\va,\vb,\vc\}$
imply that 
\[R(\alpha)R^T
(\alpha)=\brr{\cos^2
\alpha+\sin^2\alpha}\brr{\va\va^T+\vb\vb^T}+\vc\vc^T=I\]
and   hence
$\det(R(\alpha))=\pm1$.  Letting  $\lim_{\alpha\to 0}R(\alpha)=I$  and using
the continuity of the determinants, we see that $\det(R(\alpha))=1$.

We turn now to rotations in $\setR^4$. It turns out that rotations in $
\setR^4$ can be defined in a similar way to rotations in $\setR^2$ and
$\setR^3.$
We say that a linear transformation $R$ is a rotation if there exists two
dimensional subspace $\Pi$ of $\setR^4$ such that the angle between vectors
$R\vx$ and $\vx$ is a constant for all $\vx\in\Pi,$ and the angle between
vectors $R\vy$ and $\vy$ is a constant for all $\vy\in\Pi^\perp,$ the orthogonal
complement of $\Pi.$

The calculation are done in a similar
manner as we did in $\setR^3.$ Pick $\va,\vb\in\Pi$  and
$\vc,\vd\in\Pi^\perp$
such that the set $\brf{\va,\vb,\vc,\vd}$ is an orthonormal basis. Let
$R=R(\alpha,\beta)$ be the rotation
matrix with rotation's angles $\alpha$ in the plane $\Pi$ and $\beta$ in the
orthogonal complement $\Pi^\perp.$ Then
\begin{equation*}
\begin{split}
 	&R \va = \cos\alpha \va + \sin\alpha \vb, \ \
	R \vb = -\sin\alpha \va +\cos\alpha \vb,  \\
	&R \vc = \cos\beta \vc + \sin\beta \vd,  \ \
	R \vd = -\sin\beta \vc + \cos\beta \vd.
\end{split}
\end{equation*}
Set $ P = \br{\va, \vb, \vc, \vd}$ and $ Q= \br{R \va, R \vb, R \vc, R\vd}$,
since the matrix $P$ is orthogonal, $R=QP^T$ and hence
\begin{equation}
\label{eq3}
\begin{split}
R&= (R\va)\va^T + (R\vb)\vb^T + (R\vc)\vc^T + (R\vd)\vd^T \\
&=\cos\alpha\brr{\va\va^T+\vb\vb^T}+\sin\alpha\brr{\vb\va^T-\va\vb^T} \\ &+
\cos\beta\brr{\vc\vc^T+\vd\vd^T}+\sin\beta\brr{\vd\vc^T-\vc\vd^T}.
\end{split}
\end{equation}

We can now easily distinguish between two types of 4D-rotations.
If $\beta=0$, then the rotation is simple, that is, $R\vy = \vy$ for all
$\vy\in\Pi^\perp.$ Otherwise, both planes
$\Pi$ and $\Pi^\perp$ rotate simultaneously and this type is called a double
rotation.

If one doubts  whether the matrix $R$ in (\ref{eq3}) is an orthogonal
matrix
with determinant one, then the following simple calculation will convince him.
Since
\begin{equation*}
\begin{split}
 &R(\alpha,\beta)R^T(\alpha,\beta)\\
= &\left(\cos^2
\alpha+\sin^2\alpha\right)\left(\va\va^T+\vb\vb^T\right)
+\left(\cos^2\beta+\sin^2\beta\right)\left(\vc\vc^T+\vd\vd^T\right)\\=&I,
\end{split}
\end{equation*}
 $R(\alpha,\beta)$ is an orthogonal matrix and by letting $\alpha$
and $\beta$ go to zero, we get that its determinant is one.

We are now in position to extent the geometric definition of
rotations to $\setR^n$ for arbitrary positive integer $n$.
 For $n=2p$ we say
that a linear transformation $R$ is a rotation if
there exist $p$ mutual orthogonal planes $\Pi_k$ such
that the angle between the vectors $R\vx$ and $\vx$ is a constant for all
$\vx\in\Pi_k,\,k=1,...,p.$ For $n=2p+1$ we require that there are
$p$ mutual orthogonal planes $\Pi_k$ and
in addition a line $L$ is orthogonal to $\Pi_k, k=1,\ldots ,p$ such that $R$
behaves the same as in the even on the planes $\Pi_k$ and $R\vy = \vy$ for all $\vy\in
L$. The extension of formulas (\ref{eq2}) and (\ref{eq3}) to arbitrary
dimension is
obvious. In $\setR^{2p}$ there is an orthonormal basis
$\{(\va_1,\vb_1),...,(\va_n,\vb_n)\}$ such that
\begin{equation}
\label{eq11}
 R=\sum_{k=1}^p
\cos\alpha_k\brr{\va_k\va_k^T+\vb_k\vb_k^T}+\sin\alpha_k\brr{
\vb_k\va_k^T-\va_k\vb_k^T }.
\end{equation}
In odd dimension $2p+1$ there is  an orthonormal basis
$\{(\va_1,\vb_1),...,\\(\va_n,\vb_n),\vc\}$ such that
\begin{equation}
\label{eq12}
 R=\sum_{k=1}^p
\cos\alpha_k\brr{\va_k\va_k^T+\vb_k\vb_k^T}+\sin\alpha_k\brr{
\vb_k\va_k^T-\va_k\vb_k^T }+\vc\vc^T.
\end{equation}

Similarly to the rotation formulas (\ref{eq2}) and (\ref{eq3}) one can check
that (\ref{eq11}) and (\ref{eq12}) are orthogonal matrices with determinate one.

Formulas (\ref{eq11}) and (\ref{eq12}) were derived by \cite{Sch} but in a
different way.
The advantage of the derivation given here is being constructive an addition to
being  appropriate for an
elementary linear algebra course. 
Formulas (\ref{eq11}) and (\ref{eq12}) can be written in a vectors'  form
\begin{equation*}
\label{eq13}
 R\vx = \sum_{k=1}^p
\cos\alpha_k\vy_k + \sin\alpha_k \vz_k
\end{equation*}
and 
\begin{equation*}
\label{eq14}
 R\vx = \sum_{k=1}^p
\cos\alpha_k\vy_k + \sin\alpha_k \vz_k +( \vc^T\vx)\vc, 
\end{equation*}
where $\vy_k$ is the projection of vector $\vx$ on the plane $\Pi_k$ and $\vz_k$
is the rotation of $\vy_k$ by an angle $\frac{\pi}{2}$  in
the plane $\Pi_k$.

\subsection{Invariant subspaces of equiangular subspaces of rotations in
$\setR^4$}

A rotation $R$ in $\setR^4$ is called {\it equiangular rotations} or {\it
isoclinic rotations} if the planes $\Pi$ and its orthogonal complement
$\Pi^\perp$  
rotate with the same angle (see e.g. \cite{Lou, Man}).
When
$\alpha\not=\beta$, then the planes $\Pi$ and $\Pi^\perp$ are the only
invariant  subspaces under the rotation  $R$.   However,  when $\alpha=\beta$,
then there are infinitely many two dimensional invariant planes  (see e.g.
\cite{Sch}).  We shall see here that this interesting phenomenon is a simple
consequence of the formula  (\ref{eq3}) and we shall also classify all the
invariant planes. 

To see this we note that when
$\alpha=\beta$, then (\ref{eq3}) becomes
\begin{equation*}
 R=\cos \alpha I+\sin \alpha {J},
\end{equation*}
where $I$ is the identity on $\setR^4$ and
${J}=\brr{\vb\va^T-\va\vb^T+\vd\vc^T-\vc\vd^T}$. Hence a subspace $U$
of $\setR^4$ is an invariant subspace of $R$ if and only if it is invariant
subspace of the matrix ${J}$.

Now  ${J}$ is a skew-symmetric matrix satisfying  ${J}^2=-I$.
Therefore it has no real eigenvalues and this implies that any non-trivial
invariant subspace must has dimension two.
Since ${J}^2=-I$, ${\rm span}\{\vu,{J}\vu\}$ is an invariant
subspace for any $\vu\in\setR^4$. On the other hand, if $U$ is a non-trivial
invariant subspace of $R,$ then  ${J}\vu\in U$ for any $\vu\in
U$.
Hence $U$ must be spanned by these vectors.
Thus we have obtained a complete classification of the invariant subspaces of
equiangular rotations which is independent of the rotation angle $\alpha$.

It follows from formulas (\ref{eq11}) and (\ref{eq12}) that if all the angles
$\alpha_k$
 are equal, then there are infinitely many
invariant
subspaces and each one of them is spanned by a vector $\vu$ and $\sum_{k=1}^p
\brr{\vb_k\va_k^T-\va_k\vb_k^T }\vu$.

\section{Concluding Remarks}

Many vector-space textbooks use the entry-by-entry definition $c_{ij}=\sum
a_{ik}b_{kj}$ for the matrices multiplications.
The operation of multiplication of a vector $\vx$ by a matrix $A$ in
accordance (\ref{MVM}) bears in itself both geometric and algebraic properties.
Therefore a decent understanding of it should be  prior to the formal definition
of matrices multiplication.  After that the matrices' multiplication in
Cayles's spirit follows naturally. The Cayley's  definition (\ref{CD}) and the
column-row rule  (\ref{CRM}) have many
advantages. In many cases they makes the computations easier in addition to
increases the comprehension. This note emphasizes two aspects of that attitude.

\bibliographystyle{amsplain}

\end{document}